\documentclass[a4paper, 12pt, titlepage, fleqn]{article}
\usepackage{times}
\usepackage{amsthm}
\newtheorem{thm}{Theorem}

\newtheorem{lemma}[thm]{Lemma}
\newtheorem{prop}{Proposition}

\usepackage{amssymb,amsmath}

\DeclareMathOperator{\Tr}{Tr}

\begin{document}
\baselineskip=16.3pt
\parskip=14pt
\begin{center}
\section*{The Number of Rational Points On Genus 4  Hyperelliptic Supersingular Curves in Characteristic 2}

{\large 
Gary McGuire\footnote{Research supported by the Claude Shannon Institute, Science
Foundation Ireland Grant 06/MI/006} and
Alexey Zaytsev\footnote{Research supported by Science Foundation Ireland Grant 07/RFP/MATF846} 
 \\
School of Mathematical Sciences\\
University College Dublin\\
Ireland}
\end{center}

\subsection*{Abstract}

One of the big questions in the area of curves over finite fields concerns
 the distribution of the numbers of points:
Which numbers occur as the number of points on a curve of genus $g$?
The same question  can be asked of various subclasses of curves.
In this article we classify the possibilities for the number of points on 
genus 4  hyperelliptic supersingular curves over finite fields of order $2^n$, $n$ odd.
 
Keywords: curve, genus, supersingular, finite field

MSC: 14H45

\section{Introduction}
  
 Throughout this paper we let $q=2^n$, where $n$ is odd, and let
 $\mathbb{F}_q$ denote a finite field with $q$ elements.

This paper concerns the possibilities for the number of $\mathbb{F}_q$-rational points,
$N$, on  hyperelliptic supersingular   curves.
 The Serre refinement of the Hasse-Weil bound gives
 \begin{equation}\label{HWS}
 |N-(q+1)|\leq g\lfloor 2\sqrt{q} \rfloor
 \end{equation}
 which  allows a wide range of possible values for $N$.
  The typical phenomenon for supersingular curves  is that the number of points is far more 
 restricted than the general theory allows. 
 
 \newpage

To be more precise, for  curves of genus less than 4 the following results are known.

\begin{thm}
(Deuring, Waterhouse)
The  number of $\mathbb{F}_q$-rational points $N$ on a supersingular genus $1$  curve 
defined over  $\mathbb{F}_q$ satisfies
$$N-(q+1) \in \{ 0, \pm \sqrt{2q}\},$$ and all these occur.
\end{thm}

\begin{thm}
(R\"uck, Xing)
 The  number of $\mathbb{F}_q$-rational points $N$ on a 
 hyperelliptic supersingular genus $2$  curve 
defined over  $\mathbb{F}_q$ satisfies
$$N-(q+1) \in \{ 0, \pm \sqrt{2q}\},$$ and all these occur.
\end{thm}

\begin{thm}
(Oort)
There are no hyperelliptic supersingular  genus $3$ curves in 
characteristic $2$.
\end{thm}

In this paper we will prove the following theorem.

\begin{thm}\label{g4}
The  number of $\mathbb{F}_q$-rational points $N$ on a 
hyperelliptic supersingular genus $4$  curve 
defined over  $\mathbb{F}_q$ satisfies
$$N-(q+1) \in \{ 0, \pm \sqrt{2q}, \pm 2\sqrt{2q}, \pm 4\sqrt{2q} \}$$
and all these occur.
\end{thm}

Examples  show that all these
values do indeed occur, see next section.
We note that $\pm 3\sqrt{2q}$ is not a possibility.

Classifying the possible numbers of points is the same as classifying
one coefficient of the zeta function, so these results can be seen as a 
contribution towards classification of zeta functions.

Our proof uses the theory of quadratic forms in characteristic 2.
This method has previously been used in this context in van der Geer-van der Vlugt
\cite{GV}.There is aslo a discussion in Nart-Ritzenthaler \cite{NR}, see Lemma 2.2,
which restricts the number of points sufficiently for their purposes, 
but does not completely classify them.

In Section 2 we present background on curves, and in Section 3 we present
background on quadratic forms.
Section 4 presents our proof using quadratic forms.
In Section 5 we present an alternative proof, under an extra hypothesis,
using the theory of abelian varieties.
The two methods of proof are completely different.
Although the second method does not prove the full result, 
we believe it is of interest.

 \section{Curves Background}
 
 The equation 
\begin{equation}\label{sshyp4}
y^2+y= x^{9} +  ax^{5}+bx^3.
\end{equation} 
defines a hyperelliptic curve of genus 4 over $\mathbb{F}_q$,
where $a,b\in \mathbb{F}_q$.
It is shown by Scholten-Zhu \cite{SZ2} that this curve is supersingular, and
conversely, that any hyperelliptic supersingular curve of genus 4 
defined over $\mathbb{F}_q$ is isomorphic \emph{over 
the algebraic closure $\overline{\mathbb{F}_q}$}
to a curve with equation (\ref{sshyp4}).

This is not a normal form for isomorphism over $\mathbb{F}_q$.
It is shown in \cite{SZ1} (using the Deuring-Shafarevitch formula)
that any genus 4 hyperelliptic curve of 2-rank 0 defined over $\mathbb{F}_q$ has
an equation of the form
\[
y^2+y= c_9 x^{9} +c_7x^7+  c_5x^{5}+c_3x^3+c_1x.
\]
It is also shown in \cite{SZ1} that this curve is supersingular if and only if $c_7=0$.
Therefore, any hyperelliptic supersingular curve of genus 4 
defined over $\mathbb{F}_q$ is isomorphic over $\mathbb{F}_q$
to a curve with equation 
\begin{equation}\label{sshyp4a}
y^2+y= fx^{9} +  ax^{5}+bx^3+cx+d
\end{equation} 
for some constants $f,a,b,c,d \in \mathbb{F}_q$.
One needs an extension field, in general, to get an isomorphism with (\ref{sshyp4}).

For examples, when $n=11$, and $w$ is a primitive element
of $\mathbb{F}_{2^{11}}$ with minimal polynomial $x^{11}+x^2+1$,
the curve
\[
y^2+y=x^9+w^{512} x^5+w^{118} x^3
\]
has $N-(2^{11}+1)=256$, and the curve
\[
y^2+y=w^9 x^9+w^{517} x^5+w^{121} x^3+w^{24}x
\]
has $N-(2^{11}+1)=-256$.  Examples with $N-(2^{11}+1)=\pm 256$ are not common.
The curve
\[
y^2+y=x^9+w^{520} x^5+w^{117} x^3+w^{14}x
\]
has $N-(2^{11}+1)=128$
and the curve
\[
y^2+y=x^9+w^{520} x^5+w^{117} x^3+w^{15}x
\]
has $N-(2^{11}+1)=-128$.
Examples with $N-(2^{11}+1)=0$ or $\pm 64$ are plentiful.

\section{Quadratic Forms Background}\label{QF}

We now outline the basic theory of quadratic forms over $\mathbb{F}_{2}$.

Let $Q:\mathbb{F}_q \longrightarrow \mathbb{F}_{2}$ be a quadratic form.
The polarization of $Q$ is the symplectic bilinear form $B$ defined by
\[
B(x,y)=Q(x+y)-Q(x)-Q(y).
\]
By definition the radical of $B$ (denoted $W$) is 
\[
W =\{ x\in \mathbb{F}_q : B(x,y)=0 \text{  for all $y\in \mathbb{F}_q$}\}.
\]
The rank of $B$ is defined to be $n-\dim(W)$, and the first basic theorem
of this subject states that the rank must be even.

Next let $Q|_W$ denote the restriction of $Q$ to $W$, and let 
\[
W_0=\{ x\in W : Q(x)=0\}
\]
(sometimes $W_0$ is called the singular radical of $Q$).
Note that $Q|_W$ is a linear map $W\longrightarrow \mathbb{F}_2$
with kernel $W_0$.  Therefore
\begin{equation*}
\dim W_0 =
\begin{cases}
\dim(W) -1 & \text{if $Q|_W$ is onto} \\ 
\dim(W) & \text{if $Q|_W =0$ (i.e. $W=W_0$)}.%
\end{cases}%
\end{equation*}
The rank of $Q$ is defined to be $n-\dim(W_0)$.
The following theorem is well known, see \cite{GV} or \cite{LN} for example.

\begin{thm}\label{counts}
Continue the above notation.
Let $M=|\{x\in \mathbb{F}_q : Q(x)=0\}|$, and let $w=\dim(W)$.

If $Q$ has odd rank then $M=2^{n-1}$.
In this case, $\sum_{x\in \mathbb{F}_q} (-1)^{Q(x)} =0$.

If $Q$ has even rank then $M=2^{n-1}\pm 2^{(n-2+w)/2}$.
\end{thm}

\newpage

  \section{Proof of Theorem 4}
 
 Determining the value of the sum
$$S:=\sum_{x\in \mathbb{F}_q} (-1)^{\Tr(fx^{9} +  ax^{5}+bx^3+cx+d )}$$
is equivalent to determining the number of $x\in \mathbb{F}_q$ for which 
$\Tr( fx^{9} +  ax^{5}+bx^3+cx+d)=0$.
By Hilbert's Theorem 90, this is equivalent to determining the number of 
solutions in $ \mathbb{F}_q$  of   (\ref{sshyp4a}).
Indeed, if $N$ is the number of projective $\mathbb{F}_q$-rational points on 
$$y^2+y= fx^{9} +  ax^{5}+bx^3+cx+d$$
then $S=N-(q+1)$.

\begin{thm}\label{specvalues}
$S$ must take values in the set
$$\{0, \pm 2^{(n+1)/2}, \pm 2^{(n+3)/2}, \pm 2^{(n+5)/2} \}.$$
Equivalently, $N-(q+1)$ must take values in the set
$$\{0, \pm \sqrt{2q}, \pm 2\sqrt{2q},\pm  4\sqrt{2q}  \}.$$
\end{thm}

Proof:
 Squaring $S$ gives
 \[
 S^2=\sum_{x,y\in \mathbb{F}_q} (-1)^{\Tr( fx^{9} + fy^9+ ax^{5}+ay^5+bx^3+by^3+cx+cy)}.
 \]
 Substituting $y=x+u$, and for notational purposes letting $\chi(t)=(-1)^{\Tr (t)}$, we get
 \begin{eqnarray*}
 S^2&=&\sum_{x,u\in \mathbb{F}_q} \chi( fx^{9} + f(x+u)^9+ ax^{5}+a(x+u)^5+bx^3+b(x+u)^3+
 cx+c(x+u) )\\
&=& \sum_{x,u\in \mathbb{F}_q} \chi( f(x^{8}u +xu^{8}+u^{9})+  a(x^{4}u+xu^4+u^5)
+b(x^2u+xu^2+u^3)+cu)\\
&=& \sum_{u\in\mathbb{F}_q} \chi (fu^{9} +  au^{5}+bu^3+cu)\biggl(
\sum_{x\in \mathbb{F}_q} \chi(f(x^{8}u +xu^{8})+  a(x^{4}u+xu^4)+b(x^2u+xu^2) )\biggr)\\
&=& \sum_{u\in\mathbb{F}_q} \chi (fu^{9} +  au^{5}+bu^3+cu)\biggl(
\sum_{x\in \mathbb{F}_q} \chi( x^{8}[
fu +f^8 u^{64}+  a^2u^2+a^{8}u^{32}+b^{4}u^{4}+b^{8}u^{16}])\biggr).
 \end{eqnarray*}
To obtain the last line we have used the facts that
 $\chi(s+t)=\chi(s) \chi(t)$ and $\chi(t^2)=\chi(t)$.
 
 The inner sum has the form $\sum_{x\in \mathbb{F}_q} \chi( xL)$, and is a character sum over a group
 because $\chi$ is a character of the additive group of $\mathbb{F}_q$.
 This sum is therefore 0 unless $L=0$.
 Letting
 \[
 L_{f,a,b}(u)=L(u)= fu +f^8 u^{64}+  a^2u^2+a^{8}u^{32}+b^{4}u^{4}+b^{8}u^{16}
 \]
 we have
 \[
 S^2= \sum_{u\in\mathbb{F}_q} \chi (fu^{9} +  au^{5}+bu^3+cu)\biggl(
\sum_{x\in \mathbb{F}_q} \chi( x^{8}L(u))\biggr)
 \]
 and the inner sum is 0 unless $L(u)=0$.
 
 Note that $L(u)$ is a linearized polynomial, and the roots form a vector space over $\mathbb{F}_2$.
Let $W_{f,a,b}=W$ be the kernel of $L(u)$ inside $\mathbb{F}_q$, i.e.,
\[
W=\{ u\in \mathbb{F}_q : L(u)=0\}.
\]
Then $W$ is an $\mathbb{F}_2$-subspace of  $\mathbb{F}_q$ of dimension at most $6$,
because $L(u)$ has degree $64=2^6$.
 
 The next part of the proof is to observe that the dimension of $W$ must be odd.
 This is because $n-\dim (W)$ is the rank of a symplectic bilinear form, and this rank must be even.
 The form here is $B(x,y)=Q(x+y)-Q(x)-Q(y)$ where
 $Q(x)=\Tr(fx^{9} +  ax^{5}+bx^3+cx)$ is an $\mathbb{F}_2$-valued quadratic form.
 It is straightforward to check that $W$ is the radical of $B$.
 
 We now conclude that 
 $W$ is an $\mathbb{F}_2$-subspace of  $\mathbb{F}_q$ of dimension $1$ or $3$ or $5$.

 We may now write
 \[
 S^2= q \sum_{u\in W} \chi (fu^{9} +  au^{5}+bu^3+cu).
 \]
 If $Q$ is not identically 0 on $W$, then $S=0$ by Theorem \ref{counts}, 
 because $\chi$ is a non-trivial character  on $W$ and  $Q$ has odd rank.
 On the other hand, if $Q$ is identically 0 on $W$, 
 then $Q$ has even rank and by Theorem \ref{counts} we get
 \[
 S^2=q\cdot |W|.
 \]
 Because $|W|=2^w$ where $w\in \{ 1,3,5\}$ we are done.
 $\square$

\begin{section}{Jacobians of supersingular hyperelliptic curves}

In this section we prove the main theorem of this article via the theory of supersingular abelian varieties, but under an extra assumption, which most likely holds for this type of abelian variety.

Recall that a curve $C$ is called  supersingular over a finite field $\mathbb{F}_q$ if its Jacobian
${\rm Jac(C)}$ is a supersingular abelian variety. Therefore we start with some facts from the theory of abelian varieties.

Let $A$ be an abelian variety of dimension $g$ over $\mathbb{F}_q$ where $q =p^n$. The characteristic polynomial $P(X)\in \mathbb{Z}[X]$ of the Frobenius endomorphism is given by
\begin{displaymath}
P_A(X) = X^{2g} + a_{1}X^{2g-1} + \cdots+ a_gX^g +qa_{g-1}X^{g-1} \cdots + q^g.
\end{displaymath}
Recall that an abelian variety $A$ is simple if it is not isogenous to a product of abelian varieties of lower dimensions. 
In that case $P_A(X)$ is either irreducible over $\mathbb{Z}$ or $P_A(X)=h(X)^e$ where $h(X)\in \mathbb{Z}[X]$ is irreducible over $\mathbb{Z}$.

The isogeny classes are completely classified by their characteristic polynomials.
\begin{thm}(Tate) \label{Tate} 
Let $A$ and $B$ be the abelian varieties defined over $\mathbb{F}_q$. Then an abelian variety 
$A$ is $\mathbb{F}_q$-isogenous to an abelian subvariety of $B$ if and only if $P_A(X)$ divides $P_B(X)$ over $\mathbb{Q}[X]$. In particular, $P_A(X) = P_B(X)$ if and only if $A$ and $B$ are $\mathbb{F}_q$-isogenous.
\end{thm}

The property of an abelian variety being supersingular is a property of the isogeny class of this abelian variety and therefore it is completely determined by its characteristic polynomial. For example 
there is the Stichtenoth-Xing criteria \cite{Stich-Xing} that imposes conditions on the coefficients of a characteristic polynomial $P_A(X)$ of an abelian variety $A$ over a finite field $\mathbb{F}_q$ to be a supersingular abelian variety. Their result is the following.

\begin{thm}\label{St-Xing}  Let $q = p^n$ and $A$ be an abelian variety of dimension $g$ over
a finite filed $\mathbb{F}_{q}$ with the characteristic polynomial of the Frobenius endomorphism 
$$P_A(X) = X^{2g} + a_{1}X^{2g-1} + \cdots+ a_gX^g +qa_{g-1}X^{g-1} \cdots + q^g.$$
Then $A$ is supersingular if and only if $\, p^{\lceil{\frac{jn}{2}}\rceil}|a_j$ for all $1 \leq j \leq g.$ 
\end{thm}

The main object of investigation of this paper is a supersingular curve of genus $4$ 
over a finite field $\mathbb{F}_{2^n}$, where $n$ is an odd number. Therefore we start with a list of the 
characteristic polynomials of a simple supersingular abelian variety of dimension less or equal $4$ over $\mathbb{F}_{2^n}$, see \cite{KMZ}.
$$
\begin{array}{l}
X^2 \pm 2^{(n+1)/2}X+2^n\\
X^2+2^n\\
X^4\pm2^{n}X^2+2^{2n}\\
X^4\pm 2^{(n+1)/2}X^3 +2^n X^2 \pm 2^{(3n+1)/2}X +2^{2n}\\
(X^2-2^n)^2\\
X^8 \pm 2^{(n+1)/2}X^7+2^n X^6-2^{2n}X^4+2^{3n}X^2 \pm 2^{(7n+1)/2}X+2^{4n}\\
X^8+2^{4n}\\
X^8-2^{n}X^6+2^{2n}X^4-2^{3n}X^2+2^{4n}.\\
\end{array}
$$
As a direct consequence of  Tate's theorem  (Theorem  \ref{Tate}) we obtain the following proposition.
\begin{prop}\label{genus4}
Let $C$ be a supersingular curve of genus $4$ over a finite field $\mathbb{F}_{2^n}$, where $n$ is odd, such that $\#C(\mathbb{F}_{2^n})=2^n+1 \pm 3 \sqrt{2^{n+1}} $. Then 
the polynomial
$$
(X^2+2^{(n+1)/2}X+2^n)^3(X^2+2^n)
$$
is the characteristic polynomial of the Frobenius on ${\rm Jac}(C)$.
\end{prop}

We prove that under certain conditions, this type of curve is impossible.

Now, we consider the endomorphism algebra of the ${\rm Jac}(C)$, 
and consider the  element
$({\rm Frob+Ver})/2^{(n-1)/2}$
(where $Frob$ is the Frobenius endomorphism and $Ver$ is the Verschiebung).
From Theorem \ref{St-Xing} it follows that this element is algebraic:

\begin{lemma}
Let $\omega$ be a root of the characteristic polynomial of Frobenius of supersingular curve. Then the number $(\omega+\bar{\omega})/2^{(n-1)/2}$ is an algebraic number.
\end{lemma}

Now we can prove the main theorem about  the number of rational curves on the  curves
under consideration.

\begin{thm}
Let $C$ be a hyperelliptic supersingular curve of genus $4$ over a finite field $\mathbb{F}_{2^n}$, with odd $n$. If $({\rm Frob}+{\rm Ver})/2^{(n-1)/2}$ is an endomorphism of ${\rm Jac}(C)$ then
$$
\#C(\mathbb{F}_{2^n})-2^n-1 \in \{0,\, \pm 2^{(n+1)/2}\, \pm2\cdot 2^{(n+1)/2},\, \pm 4\cdot 2^{(n+1)/2}  \}.
$$
\end{thm} 

\begin{proof}
From Theorem \ref{St-Xing} it follows that ${\rm Tr(Frob_{\mathbb{F}_{2^n}})}$ is divisible by $2^{(n+1)/2}$,
hence $2^{(n+1)/2}$ divides $\#C(\mathbb{F}_{2^n})-2^n-1$. Combing this result and  the Hasse-Weil-Serre bound we get  that $\#C(\mathbb{F}_{2^n})-2^n-1=N 2^{(n+1)/2}$, with $|N| \le 4$.

Now we assume that $\#C(\mathbb{F}_q)=2^n+1\pm 3\cdot 2^{(n+1)/2}$ . Then by Proposition~\ref{genus4} and Tate's theorem it follows that
${\rm Jac(C)} \sim E_1^{3}  \times E_2$, where $E_1$ and $E_2$ are supersingular elliptic curves with  characteristic polynomials $X^2+2^{(n+1)/2}X+2^n$ and $X^2+2^n$, respectively.
The polynomial $(X-2)^3X$is the characteristic polynomial of $({\rm Frob}+ {\rm Ver})/2^{(n-1)/2}$,
and the resultant of these factors is 2. Therefore by theorem 1 in \cite{Lauter} it follows there exists a curve
$D$ such that ${\rm Jac(D)}$ is isogenous either $E_1^3$ or $E^2$ and there exists a double cover
$\sigma:C \rightarrow D$. Due to the fact the hyperelliptic involution lies in the center of an automorphism group of a hyperelliptic curve we get the following decomposition
$$
{\rm Jac}(C) \sim  {\rm Jac}(D) \times {\rm Jac}(C/ \langle \sigma \tau\rangle)
$$
and a double covering $\sigma \tau: C  \rightarrow C/ <\sigma \tau>$.
Therefore we have that either $D$ has genus $3$ or $C/ <\sigma \tau>$  has genus $3$,
which is impossible by  Hurwitz's genus formula.
\end{proof}

\end{section}


 \end{document}